\documentclass[11pt]{amsart}
\usepackage{amsmath,amsfonts,amsthm,amscd,amssymb,mathrsfs,amssymb,bm}
\usepackage{mathrsfs}
\usepackage{graphicx}
\usepackage{wrapfig}
\usepackage[all]{xy}
\newtheorem{theorem}{Theorem}

\newtheorem{proposition}[theorem]{Proposition}
\newtheorem{lemma}[theorem]{Lemma}

\newtheorem{remark}[theorem]{Remark}
\newcommand{\aaa}{\alpha}
\newcommand{\bbb}{\beta}
\newcommand{\ccc}{\gamma}
\newcommand{\ddd}{\delta}

\newcommand{\CP}{\mathbb{CP}}
\newcommand{\CC}{\mathbb{C}}

\newcommand{\RR}{\mathbb{R}}
\newcommand{\ZZ}{\mathbb{Z}}
\newcommand{\FF}{\mathbb{F}}

\newcommand{\U}{{\rm{U}}}
\newcommand{\SU}{{\rm{SU}}}
\newcommand{\LB}{{\rm{LB}}}

\newcommand{\Sing}{{\rm{Sing}\,}}

\newcommand{\ol}{\overline}

\newcommand{\lras}{\,\longrightarrow\,}

\newcommand{\set}{\,|\,}
\newcommand{\proofend}{\hfill$\square$}

\newcommand{\Aut}{{\rm{Aut}}}

\newcommand{\ms}{\mathscr}
\newcommand{\vsp}{\vspace{3mm}}

\numberwithin{equation}{section}
\numberwithin{theorem}{section}

\begin{document}
\bibliographystyle{alpha} 
\title{Deformation of LeBrun's ALE metrics with negative mass}
\author{Nobuhiro Honda}
\address{Mathematical Institute, Tohoku University,
Sendai, Miyagi, Japan}
\email{honda@math.tohoku.ac.jp}
\thanks{The author was partially supported by the Grant-in-Aid for Young Scientists  (B), The Ministry of Education, Culture, Sports, Science and Technology, Japan. }
\begin{abstract}
In this article 
we investigate
deformations of a scalar-flat K\"ahler metric on the total space of  
complex line bundles over $\CP^1$ constructed by C.\,LeBrun.
In particular, we find that the metric is included in a one-dimensional
family of such metrics on the four-manifold, where the complex structure in the deformation is  not the standard one.
\end{abstract}
\maketitle
\setcounter{tocdepth}{1}
\vspace{-5mm}


\section{Introduction}
In 1988 C.\,LeBrun \cite{LB88} explicitly constructed
an example of anti-self-dual (ASD) K\"ahler metric on the total space of the complex line bundle $\ms O(-n)$
over $\CP^1$, which is asymptotically 
locally Euclidean (ALE) and whose mass is negative when $n>2$.
Significance of the metric is not only 
in that it provides counter-example to 
the generalized positive action conjecture, but also
in that, it naturally appears 
in a typical example 
(\cite[Section 5]{LB91}) of degeneration
of compact ASD manifolds
 as one of the two pieces
(see also \cite{DF89,KS01}).
In this degeneration, the other piece is an ALE hyper-K\"ahler manifold
constructed by Gibbons-Hawking \cite{GH78} and
Hitchin-Kronheimer \cite{H79,Kr89a}.

Because the LeBrun metrics can be thought as a natural generalization
of the Burns metric and the Eguchi-Hanson metric on $\ms O(-1)$ and $\ms O(-2)$
respectively, 
and since these two metrics are rigid as ALE ASD metrics, 
one might think that  the LeBrun's metrics could not be
deformed as an ASD structure.
However, in a very recent work, 
by establishing an index theorem for the deformation complex on compact
ASD orbifolds, 
J.\,Viaclovsky \cite{V} 
has shown that the versal family of the LeBrun's ALE ASD structure
on $\ms O(-n)$ is non-trivial, and 
that the moduli space of ASD structures near 
the LeBrun's one is at least $(4n-12)$-dimensional.
The purpose of the present paper is to 
answer some questions which naturally arise
from that work.

We recall that from the ALE condition,
LeBrun's metric can be conformally compactified
by adding one point at infinity,
and consequently we obtain an ASD structure on a compact  orbifold 
$\widehat{\ms O(-n)}$.
We call this ASD orbifold as the {\em LeBrun orbifold}.
In Section \ref{s:Kura} by making use of
 the twistor space of the LeBrun orbifold,
we reprove that the parameter space
of the versal family for
the LeBrun orbifold is 
smooth and real $(4n-8)$-dimensional.
Here we are considering versal family of ASD structures
on the fixed orbifold $\widehat{\ms O(-n)}$.
Because the LeBrun orbifold has an effective 
U(2)-action \cite{LB88}, the parameter space of the
above versal family also has a natural
U(2)-action.
Then by determining this U(2)-action and classifying all U(2)-orbits whose dimension is less than four,
we prove the following result:

\begin{theorem}\label{thm:smry1}
Let $n\ge 3$, $B$ be an open neighborhood of 
the origin in $\RR^{4n-8}$,
 and $\{[g_t]\set t\in B\subset\RR^{4n-8}\}$ 
be
the versal family of ASD structures
on the orbifold $\widehat{\ms O(-n)}$
for the LeBrun's ASD structure, so that $[g_0]$ is
equal to the LeBrun's ASD structure.
If we take $B$ sufficiently small, the following holds.
\begin{enumerate}
\item[\em (i)] If $n=3$, for any $t\in B$ with $t\neq 0$,
we have $\Aut_0[g_t]\simeq\U(1)$,
and the moduli space is 1-dimensional at 
the point $[g_t]$.
\item[\em (ii)] If $n\ge 3$, 
there exist  $\U(2)$-invariant, mutually 
disjoint connected subsets $B_1,\cdots,B_{[n/2]}$ of $B$ such that $t\in B_1\cup\cdots\cup B_{[n/2]}$ implies
$\Aut_0[g_t]\simeq \U(1)$. 
Moreover, the moduli space of these
$\U(1)$-invariant 
ASD structures is $1$-dimensional at $[g_t]$
if $t\in B_1$,  $3$-dimensional
at $[g_t]$ if $t\in B_2\cup B_3\cup\cdots\cup B_{[n/2]}$.
\item[\em (iii)] 
If $n>3$ and $t\not\in B_1\cup B_2\cup\cdots\cup B_{[n/2]}$, then $\Aut_0[g_t]=\{e\}$ and the dimension of the moduli space is $(4n-12)$-dimensional at $[g_t]$.
\item[\em (iv)]
If $n=4$, there exists another $U(2)$-invariant
connected subset $B_0$ of $B$ for which the following holds:
$t\in B_0$ implies
$\Aut_0[g_t]\simeq \SU(2)$,
and the moduli space of these
$\SU(2)$-invariant ASD structures is 
1-dimensional at $[g_t]$.
Further if $t\not\in B_0\cup B_1\cup B_2$,
then the moduli space is $4$-dimensional at $[g_t]$.
\end{enumerate} 
\end{theorem}
\noindent
Here, for a conformal structure $[g]$,
$\Aut_0[g]$ denotes the identity component of 
the conformal automorphism group of $[g]$,
and
for a number $k$, $[k]$ means the 
largest integer not greater than $k$.
Theorem \ref{thm:smry1} classifies
all small deformations of the LeBrun's ASD structure which are equivariant with respect to
a subgroup of $\U(2)$ of positive dimension,
and
gives an answer to a question by Viaclovsky \cite[1.4.\,Question (3)(4)]{V}
regarding deformations of the LeBrun orbifold.
For the $\SU(2)$-equivariant deformation of the metrics on $4\CP^2$ obtained in (iv)
of the theorem, 
it would be interesting to find concrete 
description of them, under the work of Hitchin \cite{H95}
concerning $\SU(2)$-invariant ASD metrics in general.

In Section \ref{s:sfk} we study deformations
of LeBrun's ASD orbifold which preserves
K\"ahlerity on the smooth locus, again by using twistor space.
The key for such investigation is a theorem 
of Pontecorvo \cite{Pont92} which expresses
the K\"ahlerity of an ASD structure in terms of
certain divisor on the twistor space.
Especially we prove the following result:

\begin{theorem}\label{thm:sfk}
For any $n\ge 3$, 
on the 4-manifold $\ms O(-n)$,
there exists a one-dimensional smooth family 
$\{(J_t,g_t)\}$
of  complex structures and ALE, ASD K\"ahler metrics, which satisfies the 
following properties:
\begin{enumerate}
\item[\em (i)] $g_0$ coincides with the LeBrun metric, and $J_0$ is the standard complex structure,
\item[\em (ii)] if $t\neq 0$, $g_t$ is not 
conformal to the LeBrun metric.
Further 
 the complex surface $(\ms O(-n),J_t)$ is
biholomorphic to an affine surface in $\CC^{n+1}$.
Furthermore, the K\"ahler surface
$(\ms O(-n),J_t,g_t)$ admits
a non-trivial $\U(1)$-action.
\end{enumerate} 
\end{theorem}
\noindent
In particular  if $t\neq 0$ 
the complex structure $J_t$ on $\ms O(-n)$
is different from the standard one.
More explicitly, the affine surface in $\CC^{n+1}$ in the theorem can be concretely obtained as follows.
Let $\FF_{n-2}:=\mathbb P(\ms O(n-2)\oplus\ms O)$ be the ruled surface over $\CP^1$,  and $\Gamma$ and $h$ the unique negative section
and a fiber of the ruling respectively.
Then the linear system $|\Gamma + (n-1)h|$ induces an embedding $\FF_{n-2}
\subset\CP^{n+1}$. 
Thus if we remove a generic hyperplane section (which is 
a $(+n)$-rational curve) from the image of $\FF_{n-2}$, we get an 
affine surface in $\CC^{n+1}$.
This is nothing but the affine surface in the theorem.
Note that if $p:\CC^{n+1}\to \CP^n$ denotes the projection from the origin,
the total space $\ms O(-n)$ can be realized as the minimal resolution of 
the cone for the projection $p$ over the rational normal curve in $\CP^n$.
Then the smooth affine surface discussed above is obtained by varying the 
defining (quadratic) equations of the cone and taking a simultaneous resolution
of the cone singularity in the family (see \cite[\S 8]{Pin74}) .

%
%
%
%

Theorem \ref{thm:sfk} provides a partial answer to 
a Viaclovsky's question \cite[1.4.\,Question (2)]{V} concerning scalar-flat K\"ahler deformations of the LeBrun metric.
From the proof the family of the K\"ahler metrics
in Theorem \ref{thm:sfk}
is obtained from the U(1)-equivariant family
over the subset $B_1$ in Theorem \ref{thm:smry1}
by restricting onto any 1-dimensional linear subspace in $B_1$.
We also investigate other non-trivial
U(1)-equivariant deformations of the
LeBrun's orbifold, 
which are over $B_2,B_3,\cdots, B_{[n/2]}$ in Theorem \ref{thm:smry1},
and 
show that they do not preserve K\"ahlerity
of the metric, in contrast with the above one.
Further we also observe that
for these deformations the corresponding 
twistor spaces are non-Moishezon.
This would be  natural in light of similar
phenomena in the case of twistor spaces on $n\CP^2$. 

After writing this paper, 
Michael Lock and Jeff Viaclovsky \cite{MV2012} extended the index theorem in \cite{V} to  general compact ASD orbifolds with cyclic quotient
singularities, and showed for example that the ALE SFK metrics
constructed by Calderbank-Singer
\cite{CS2004}
on the minimal resolution of the quotient $\CC^2/\Gamma$,
$\Gamma$ being a cyclic group,  admit a non-trivial deformation as ALE ASD metrics.
But it is not straightforward to see that the method  used in this paper can be applied to the twistor spaces of their spaces,
since the singularities on the twistor spaces are not Gorenstein any more (i.e.\,the canonical divisor is not a Cartier divisor),
which makes the key divisor $S$ non-Cartier.

\vsp
\noindent{\bf Notation.}
We write $\FF_n$ for the ruled surface $\mathbb P(\ms O(n)\oplus\ms O)$ over $\CP^1$.
If $X$ is a subset of a twistor space,
we denote by $\ol X$ for the image of $X$
under the real structure.

\vsp
\noindent
{\small
{\bf {Acknowledgement}}
I would like to thank Jeff Viaclovsky for 
letting me know his latest result concerning 
deformations of the LeBrun's ALE metric,
and also for many helpful suggestions.
We also thanks Kazuo Akutagawa,
Akira Fujiki,  Henrik Pedersen, Yat Sun Poon, and Carl Tipler
for  useful comments and suggestions which substantially improve
this article.
%
%
}

\section{The Kuranishi family of the twistor space}\label{s:Kura}
\subsection{The twistor space of the LeBrun orbifold}\label{ss:recall}
First we briefly recall the LeBrun's ALE metric
with negative mass from \cite{LB88}.
For more details, one can also consult a paper by Viaclovsky \cite[Sections 2.3, 5.2]{V10}.
Fix any integer $n>2$.
On $\CC^2$, the metric is written as
\begin{align}\label{gLB}
g_{\,\LB}:= \frac{dr^2}
{\left(
1-\frac 1{r^2}
\right)
\left(
1+\frac{n-1}{r^2}
\right)}
+r^2
\left[
\sigma_1^2 + \sigma_2^2 +
\left(
1-\frac 1{r^2}
\right)
\left(
1+\frac{n-1}{r^2}
\right)
\sigma_3^2
\right],
\end{align}
where $r$ is the Euclidean distance form the origin,
and $\sigma_1,\sigma_2,\sigma_3$ are  left-invariant
coframe of  $\SU(2)=S^3$.
Clearly this metric has singularities at the unit sphere.
Let $\zeta:=e^{2\pi i/n}$ and consider the action of the cyclic group $\ZZ_n=\ZZ/n\ZZ$
on $\CC^2$ generated by 
\begin{align}\label{act01}
(z,w)\longmapsto (\zeta z,\zeta w).
\end{align}
After dividing $\CC^2$ by this action
and resolving the resulting singularity at the origin, the metric
\eqref{gLB} defines a non-singular Riemannian metric
on the total space of the holomorphic line bundle
$\ms O(-n)\to \CP^1$, which is locally asymptotically Euclidean (ALE) at infinity.
Moreover, the metric is K\"ahler with respect to
the complex structure on $\ms O(-n)$.
In particular, it is anti-self-dual (ASD).
Furthermore from the ALE property,
after an appropriate conformal change,
the  metric extends to a one-point
 compactification
$\widehat{\ms O(-n)}=\ms O(-n)\cup\{\infty\}$
 as an orbifold ASD structure.
For brevity we call the orbifold
$\widehat{\ms O(-n)}$ equipped with
this ASD structure 
as the {\em LeBrun orbifold}.

We also recall that the isometry group of the LeBrun orbifold is the unitary group U(2), and 
the U(2)-action on the open subset $\ms O(-n)$ is realized  from 
the natural U(2)-action on $\CC^2$ 
through the quotient by $\ZZ_n$ and the minimal resolution.

The twistor space of the LeBrun orbifold is implicitly constructed in 
his different paper \cite[Section 3]{LB91}, and we now recall 
the construction,  according to \cite{HonCMP}.
Let $n\ge 3$ be an integer as above.
Over $\CP^1\times\CP^1$, take a rank-3 vector bundle 
$
\ms E_n:= \ms O(n-1,1)\oplus\ms O(1,n-1) \oplus \ms O,
$
and consider the associated $\CP^2$-bundle $\mathbb P(\ms E_n)\to\CP^1\times\CP^1$.
Let $X\subset \mathbb P(\ms E_n)$ be a hypersurface defined by
\begin{align}\label{def:X}
xy = (u-v)^n t^2,
\end{align}
where $(u,v)$ are  non-homogeneous  coordinates on $\CP^1\times\CP^1$,
and $(x,y,t)$ are fiber coordinates on the bundle $\ms E_n$.
This is an equation which takes values  in the line bundle $\ms O(n,n)$ over $\CP^1\times\CP^1$.
$X$ is equipped with an anti-holomorphic 
involution. See \cite{HonCMP} for its 
concrete form.
Next define  divisors $E$ and $\ol E$, and a curve $L'_{\infty}$ lying on $\mathbb P(\ms E_n)$ by
\begin{align}\label{def:el}
E=\{x=t=0\},\quad \ol E= \{y=t=0\}
\quad{\text{and}}\quad L'_{\infty}=\{x=y=u-v=0\}.
\end{align}
These are included in $X$, and $E$ and $\ol E$ are sections
of the projection $X\to \CP^1\times\CP^1$.
$L'_{\infty}$ is non-singular and isomorphic to
$\CP^1$.
Moreover we have $E\cap \ol E=\emptyset$ and $(E\cup \ol E)\cap L'_{\infty}=\emptyset$.
The threefold $X$ has $A_{n-1}$-singularities along the curve $L'_{\infty}$,
and this is exactly the singular locus of $X$.
By looking the normal bundle in $X$, the  section $E$ 
can be blown-down in a unique way to $\CP^1$ along a projection $\CP^1\times\CP^1\to \CP^1$, and the same for $\ol E$.
Let $\mu:X\to Z_{\LB}$ be the blowdown of $E\cup\ol E$ obtained this way, and put $C=\mu(E)$ and $\ol C=\mu(\ol E)$ for the image rational curves.
We have $N_{C/Z_{\LB}}\simeq N_{\ol C/Z_{\LB}}\simeq
\ms O(1-n)^{\oplus 2}$ for the normal bundles.
These curves play significant role for studying deformations
of $Z_{\LB}$.
We write the curve $\mu(L'_{\infty})$ by $L_{\infty}$. 
Then the variety $Z_{\LB}$ is exactly the twistor space
of the LeBrun orbifold $\widehat{\ms O(-n)}$,
 the curve $L_{\infty}\subset Z_{\LB}$ 
is  the twistor line over the orbifold point $\infty$,
and the anti-holomorphic involution 
on $Z_{\LB}$ induced from that on $X$ is
the real structure 
(see \cite[Theorem 3.3]{HonCMP}).
For a later purpose we further define other two divisors on $X$ as
\begin{align}\label{D}
D'=\{x=u-v=0\}, \quad 
\ol D'=\{y=u-v=0\}.
\end{align}
These are over the diagonal $\Delta:=\{u=v\}\subset\CP^1\times\CP^1$,
and we have $D'\cap\ol D'=L'_{\infty}$.
The two divisors $D'$ and $\ol D'$ are non-singular
and biholomorphic to the ruled surface $\FF_n$.
Note that these are not Cartier divisors.
Because $\Delta$ is a $(1,1)$-curve,
the blowdown $\mu$ induces a biholomorphic
map $D'\to\mu(D')$ and $\ol D'\to\mu(\ol D')$.
We write the images by $D=\mu(D')$ and $\ol D=\mu(\ol D')$.
These are again non-Cartier divisors on $Z_{\LB}$.
Since $X$ is a hypersurface in a smooth space, the canonical line bundle
$K_X$ naturally makes sense by adjunction formula, and
we obtain 
\begin{align}\label{cbl}
K_X &\simeq \pi^*\ms O(-2,-2)  - (E+\ol E),
\end{align}
where $\pi:X\to\CP^1\times\CP^1$ denotes the
projection.
Noting $D'+\ol D'\in |\mu^*\ms O(D+\ol D) 
- E - \ol E|$, 
this implies $D + \ol D\in |K^{-1/2}|$
on $Z_{\LB}$.
By Pontecorvo's theorem \cite{Pont92},
this divisor gives a reason for the LeBrun metric to be K\"ahler, in terms of  twistor space.

The U(2)-action on the LeBrun orbifold naturally induces
a holomorphic U(2)-action on the twistor space $Z_{\LB}$, which clearly preserves the twistor line $L_{\infty}$.
Looking at the normal bundles, we readily see that 
the divisors $D, \ol D$ and the curves $C$ and $\ol C$
on $Z_{\LB}$
are invariant under this U(2)-action.

\subsection{Locally trivial deformations of the twistor space}
\label{ss:dfm}
Before going to the actual computations for the twistor spaces, we briefly recall
well-known facts regarding deformation theory for
general compact complex varieties.
For a complex variety $Y$ which may have singularities, 
let $\Omega^1_Y$ be the sheaf of K\"ahler differentials on $Y$
as usual,
and we define the tangent sheaf of $Y$ as
$$\Theta_Y:=\ms H om_{\ms O_Y} (\Omega^1_Y,\ms O_Y).$$
If $Y$ is a hypersurface in a smooth space $V$
(just as our $X$ in $\mathbb P(\ms E_n)$), 
in terms of local defining equation $f=0$ of $Y$ in $V$,
this can be concretely written as
\begin{align}\label{tansh}
\Theta_Y= \left\{v|_Y\set v\in \Theta_V, \, v(f)=0\right\}.
\end{align}
Then if $Y$ is compact, the Zariski tangent space of the Kuranishi family of {\em locally trivial} deformations of $Y$ is
identified with the cohomology group $H^1(\Theta_Y)$, and the obstruction space is
$H^2(\Theta_Y)$.
(In this article we do not need to consider general deformations
which are not locally trivial, and so we do not
need Ext-groups.)
In particular, if $H^2(\Theta_Y)=0$, the parameter space
of the Kuranishi family is identified with an open neighborhood
of the origin in $H^1(\Theta_Y)$.
 
For the present twistor space, we have the following
\begin{proposition}\label{prop:unobs}
For the twistor space $Z_{\LB}$ of the LeBrun orbifold $\widehat{
\ms O(-n)}$, we have 
$$H^2(\Theta_{Z_{\LB}})=0.$$
Moreover in terms of the $\U(2)$-invariant rational curves
$C$ and $\ol C$ in $Z_{\LB}$, 
we have a $\U(2)$-equivariant isomorphism
\begin{align}\label{isom1}
H^1(\Theta_{Z_{\LB}}) \simeq H^1 (N_{C/Z_{\LB}}\oplus
N_{\ol C/Z_{\LB}}).
\end{align}
In particular, $h^1(\Theta_{Z_{\LB}})=4n-8$
(since $N_{C/Z_{\LB}}\simeq N_{\ol C/Z_{\LB}}\simeq\ms O(1-n)^{\oplus 2}$). 
\end{proposition}

\proof
We imitate the calculations given in \cite[Section 1.2]{Hon07-1}.
In this proof for simplicity we write $Z$ for
$Z_{\LB}$.
Let $\Theta_{X,\,E+\ol E}$ denote
the subsheaf of the tangent sheaf $\Theta_X$ consisting
of vector fields which are tangent to the divisor $E\cup\ol E$.
(Since $X$ is smooth at $E\cup \ol E$, this naturally makes sense.)
We define the subsheaf  $\Theta_{Z,\,C+\ol C}$ of 
$\Theta_Z$ in a similar way. 
Then a computation using local coordinates shows that
the blowdown $\mu:X\to Z$ induces a
U(2)-equivariant isomorphism
$
\Theta_{X,\,E+\ol E}\simeq\mu^*\Theta_{Z,\,C+\ol C}.
$
This induces an equivariant isomorphism
\begin{align}\label{isom2}
H^i(\Theta_{X,\,E+\ol E})\simeq
H^i(\Theta_{Z,\,C+\ol C}),\quad i\ge 0.
\end{align} 
From the normal bundles of $E$ and $\ol E$ in $X$,
we readily obtain $H^i(\Theta_{X,\,E+\ol E})\simeq
H^i(\Theta_X)$ for any $i\ge 0$.
Hence from \eqref{isom2} we have an equivariant isomorphism
\begin{align}\label{isom3}
H^i(\Theta_{Z,\,C+\ol C})\simeq H^i(\Theta_X),
\quad i\ge 0.
\end{align}
For computing the RHS, let $\pi:X\to \CP^1\times \CP^1=:Q$ be the projection as in \eqref{cbl}, and consider the natural homomorphism
 $d\pi:\Theta_X\to \pi^*\Theta_Q$.
Since $\pi$ is clearly submersion outside 
the line $L_{\infty}$, the support of the cokernel sheaf
for $d\pi$ is contained in $L_{\infty}$.
Further for a point $x\in L_{\infty}$,
the image of the differential $(d\pi)_x:
T_xX\to T_{\pi(x)} Q$ is readily seen to be
the subspace $T_{\pi(x)}\Delta$, where $\Delta$ is the diagonal
of $Q$ as before.
Therefore the cokernel sheaf of $d\pi$ is
exactly the normal sheaf $N_{\Delta/Q}$ under
the identification $L_{\infty}\simeq \Delta$ by $\pi$, and if we write $\ms F$ for the image
sheaf of the homomorphism $d\pi$, we obtain an exact sequence
\begin{align}\label{ses0}
0\lras\ms F\lras \pi^*\Theta_Q \lras N_{\Delta/Q} \lras 0.
\end{align}
As the map $H^0(\pi^*\Theta_Q)\to H^0(N_{\Delta/Q})$
is clearly surjective, 
this sequence easily implies 
\begin{align}\label{van1}
H^i(\ms F) = 0,\quad i>0.
\end{align}
On the other hand for the kernel sheaf $\Theta_{X/Q}$ of $d\pi$,
which consists of  vertical vector fields,
noting $\Theta_{X/Q}\simeq \ms
O_X(E+\ol E)$ from an obvious vertical vector field,
and taking the direct image $\pi_*$ of
the standard exact sequence
$0\lras \ms O_X \lras \ms O_X(E+\ol E) \lras
\ms O_E(E) \oplus \ms O_{\ol E}(\ol E) \lras 0$,
we obtain $H^i(\Theta_{X/Q})=0$ for $i>0$.
Hence from the exact sequence
\begin{align}\label{ses:v}
0 \lras \Theta_{X/Q} \lras \Theta_X \lras \ms F
\lras 0
\end{align}
we obtain 
$
H^i(\Theta_{X})\simeq H^i(\ms F)$ for any  $i>0.
$
Hence by \eqref{van1} and \eqref{isom3} we get $H^i(\Theta_X)= H^i(\Theta_{Z,\,C+\ol C})=0$
for any $i>0$.
Therefore the standard exact sequence 
\begin{align}\label{ses1}
0 \lras \Theta_{Z,\,C+\ol C} \lras
\Theta_Z \lras N_{C/Z}\oplus N_{\ol C/Z} 
\lras 0
\end{align}
induces the required isomorphism
$H^1(\Theta_Z)\simeq H^1(N_{C/Z}\oplus N_{\ol C/Z} )$ as well as the vanishing $H^2(\Theta_{Z}) =0$.
The last isomorphism is clearly U(2)-equivariant,
since all isomorphisms and the exact sequences
we have used are clearly U(2)-invariant.
\proofend

\begin{remark}{\em
For the LeBrun twistor space on $n\CP^2$ constructed in \cite{LB91},
there exist similar curves
$C_0$ and $\ol C_0$ 
and 
the cohomology group $H^1(\Theta_Z)$ is 
a direct sum of $H^1(N_{C_0/Z}\oplus N_{\ol C_0/Z})$ 
with another cohomology group (see the exact sequence (1.14) in \cite{Hon07-1}). 
The latter cohomology group precisely corresponds
to deformations as LeBrun twistor spaces.
In the present case this cohomology vanishes
as in the above proof,
and all non-trivial deformations yield non-LeBrun
twistor spaces.
}
\end{remark}

Proposition \ref{prop:unobs} means that 
the parameter space of the Kuranishi family 
of locally trivial deformations of $Z_{\LB}$
may be identified with a neighborhood of
the origin in the cohomology group $H^1 (N_{C/Z_{\LB}}\oplus
N_{\ol C/Z_{\LB}})$, which is $(4n-8)$-dimensional
over $\CC$.
Deformations as twistor spaces can be obtained
by restricting the Kuranishi family to the real locus
$H^1 (N_{C/Z_{\LB}}\oplus
N_{\ol C/Z_{\LB}})^{\sigma}$
in the neighborhood,
where $\sigma$ denotes the real structure
of $Z_{\LB}$.
We call this restricted family as {\em the versal family of  twistor spaces} for $Z_{\LB}$, and the corresponding family of
ASD structures
on $\widehat{\ms O(-n)}$ as {\em the versal family of ASD structures} for the LeBrun orbifold.
We note that there is a natural U(2)-equivariant isomorphism
\begin{align}
H^1(N_{C/Z_{\LB}})\simeq
H^1 (N_{C/Z_{\LB}}\oplus
N_{\ol C/Z_{\LB}})^{\sigma},
\end{align} 
which sends an element $\eta\in H^1(N_{C/Z_{\LB}})$ to the pair $(\eta,\,\ol\sigma^*\eta)$ in the real diagonal.
Therefore as far as we are concerned with  the versal family of  twistor spaces or 
ASD structures,
the U(2)-action on one half $$H^1(N_{C/Z_{\LB}})\simeq
\CC^{2n-4}\simeq \RR^{4n-8}$$ is
fundamental, which we next discuss.

\subsection{Explicit form of the $\U(2)$-action on $H^1$.}\label{ss:U2}
For expressing the result we consider the natural representation space $\CC^2$ 
of U(2) (acted by the multiplication of matrices), and 
for each non-negative integer $m$ we write $S^m\CC^2$ for the
$m$-th symmetric product,
where $S^0\CC^2$ means the trivial representation on $\CC$.
For convenience we promise $S^m\CC^2=0$ if $m<0$.
Let $\CC_l$ be a 1-dimensional representation
of U(2) obtained by
multiplying the $l$-th power of the determinant, and 
we write $$S^m_l\CC^2:= S^m\CC^2\otimes_{\CC}\CC_l.$$
Of course we have
$
\dim_{\CC} S^m_l\CC^2 = m+1
$
for any $m\ge 0$.
This is an  irreducible representation of $\U(2)$ for any $m,l\ge 0$.
Under these notations we have
\begin{proposition}\label{U2onH1}
Suppose $n\ge 3$.
Then under the above notation, the {\rm U(2)}-action on $H^1 (N_{C/Z_{\LB}})$ is
equivalent to the direct sum
\begin{align}\label{U2onH12}
S^{n-2}_1 \CC^2 \oplus S^{n-4}_2 \CC^2.
\end{align}
(Note that the second direct summand vanishes
when $n=3$.)
\end{proposition}

For the proof of Proposition \ref{U2onH1}
 we first recall that the U(2)-action on the
open subset $\ms O(-n)$ of the LeBrun orbifold
is induced from the natural U(2)-action 
on $\CC^2$ via the quotient and minimal resolution.
From the $\ZZ_n$-action on $\CC^2$ in \eqref{act01}, 
we can use the power $z^n=:\xi$   as a fiber coordinate of the line bundle $\ms O(-n)$ on an affine open subset of $\CP^1$.
If we put $u:= w/z$, which is a coordinate
on  the affine subset,
the pair $(\xi,u)$ can be used as coordinates on 
an open subset of $\ms O(-n)$.
Under these coordinates
the U(2)-action on $\ms O(-n)$  
is explicitly given as
\begin{equation}\label{act1}
(\xi,\,u)
\stackrel{A}\longmapsto
\left((\aaa+\bbb u)^n \xi,\,\,
\frac{\gamma+\delta u}{\aaa + \bbb u}\right),
\quad
A=\begin{pmatrix}
\aaa & \bbb \\
\ccc & \ddd
\end{pmatrix}\in \U(2).
\end{equation}

\vsp
\noindent
{\em Proof of Proposition \ref{U2onH1}.}
We again write $Z$ for $Z_{\LB}$. 
We first note that via the twistor fibration,
the divisor $D$ minus the  line $L_{\infty}$
can be U(2)-equivariantly identified with the open subset $\ms O(-n)$,
while the curves $C$ and $\ol C$
are identified with the zero-section of 
the line bundle $\ms O(-n)$.
From the inclusions $C\subset
D\subset Z$, we have 
the standard exact sequence 
\begin{align}\label{ses100}
0 \lras N_{C/D} \lras N_{C/Z} \lras
N_{D/Z}|_C \lras 0
\end{align}
for the normal bundles, which is U(2)-equivariant.
Further since $D+\ol D\in |K_Z^{-1/2}|$
we have
$$
K_D\simeq K_Z + D|_D
\simeq (-2D-2\ol D +D)|_D\simeq (-D-2\ol D)|_D
$$
and hence, since $\ol D\cap C=\emptyset$, 
by restricting this to $C$,
$[D]|_C\simeq K_D^{-1}|_C$.
Therefore 
\begin{align}
N_{D/Z}|_C \simeq [D]|_C\simeq K_D^{-1}|_C
\simeq N_{C/D}\otimes K^{-1}_C,
\end{align}
where the last isomorphism is from adjunction formula.
All these isomorphisms are clearly U(2)-equivariant.
Therefore since $N_{C/D}\simeq\ms O(-n)$
and $N_{C/D}\otimes K^{-1}_C\simeq\ms O(2-n)$,
as $n>2$, from \eqref{ses100} we get an equivariant exact sequence
\begin{align}
0 \lras H^1(N_{C/D}) \lras H^1(N_{C/Z}) \lras
H^1(N_{C/D}\otimes K^{-1}_C) \lras 0.
\end{align}
(This is not true if $n=2$.)
From this 
we obtain a U(2)-equivariant isomorphism
\begin{align}\label{H1-1}
H^1(N_{C/Z})
\simeq H^1(N_{C/D}) \oplus H^1(N_{C/D}\otimes K^{-1}_C). 
\end{align}
Since the standard open covering of $C=\CP^1$ is not U(2)-invariant, it seems difficult to compute the
action on $H^1$ by using \v{C}ech  cohomology
(as we did in \cite{Hon07-1}).
So we convert it to that on $H^0$ by Serre duality.
Namely from \eqref{H1-1}, the U(2)-action on 
$H^1 (N_{C/Z})$ can be identified with 
the dual of the U(2)-action on 
$H^0(N^{-1}_{C/D}\otimes K_C)\oplus H^0(N^{-1}_{C/D}\otimes K^2_C)$.
As $\U(2)\subset {\rm O}(4)$,
the dual action is equivalent to the original one.
So we compute the U(2)-action on 
the two direct summands.
For these we use the above coordinates $(\xi,u)$.

We first compute the U(2)-action on the former space
$H^0(N^{-1}_{C/D}\otimes K_C)$.
We put $U:=\CP^1\backslash\{(0:1)\}$, where the coordinate
$u$ is valid.
We use the 1-form $d\xi$ as a frame of 
the co-normal bundle $N_{C/D}^{-1}$ over $U$.
 For $A\in U(2)$,
 by \eqref{act1}, we have
\begin{align}
A^*d\xi &= d((\aaa + \bbb u)^n \xi) \notag\\
&= n (\aaa + \bbb u) ^{n-1} \bbb \xi\,du + 
(\aaa + \bbb u) ^{n} d\xi.
\end{align}
So over  the zero-section $\{\xi=0\}$,
we have 
\begin{align}\label{act2}
A^*d\xi = (\aaa + \bbb u)^n d\xi .
\end{align}
On the other hand, for the frame of the canonical bundle
$K_C$ on $U$, we use the 1-form $du$.
For the pull-back of $du$ under $A$, we have
\begin{align}\label{act3}
A^*du = d
\left(
\frac{\gamma+\delta u}{\aaa + \bbb u}
\right)
= \frac{\aaa\ddd-\bbb\ccc}{(\aaa + \bbb u)^2} du.
\end{align} 
By \eqref{act2} and \eqref{act3}, the U(2)-action on 
the line bundle $N_{C/D}^{-1}\otimes K_C$ is given by
\begin{align}\label{act4}
d\xi\otimes du
\stackrel{A^*}\longmapsto 
(\aaa\ddd-\bbb\ccc)\,(\aaa + \bbb u)^{n-2} 
\,d\xi\otimes du.
\end{align}  
Since $\deg (N^{-1}_{C/D}\otimes K_C )=n-2$,
any global section of $N^{-1}_{C/D}\otimes K_C$ can be written as
$P(u)\, d\xi \otimes du$ for some polynomial $P(u)$ with
$\deg P(u)\le n-2$.
For this section, by \eqref{act4}, we obtain
\begin{align}\label{key2}
P(u)\,d\xi \otimes du
\stackrel{A^*}\longmapsto 
(\aaa\ddd-\bbb\ccc)\,
\left\{(\aaa + \bbb u)^{n-2}\,P\left(
\frac{\gamma+\delta u}{\aaa + \bbb u}
\right)
\right\} 
d\xi \otimes du.
\end{align}
The ingredient of the brace is a polynomial whose degree is at most $ (n-2)$,
and
the assignment 
$P(u)\mapsto (\aaa + \bbb u)^{n-2}P((\ccc + \ddd u)/(\aaa + \bbb u) ) 
$
is exactly the $(n-2)$-th symmetric product of the natural representation
of U(2).
Therefore noting the determinant in \eqref{key2},
the U(2)-action on $H^0(N^{-1}_{C/D}\otimes K_C)$
is equivalent to $S^{n-2}\CC^2\otimes \CC_1=S_1^{n-2}\CC^2$.
Thus we obtain the first direct summand in \eqref{U2onH12}.

The  U(2)-action on the latter space
$H^0(N^{-1}_{C/D}\otimes K^2_C)$
can be readily obtained from the above computations if we notice that 
$d\xi \otimes (du)^2$ can be used as a frame over $U$,
instead of $d\xi\otimes du$.
Namely, by taking the tensor product of 
\eqref{act2} with the square of \eqref{act3}, 
the U(2)-action on $H^0(N^{-1}_{C/D}\otimes K^2_C)$
is exactly $S^{n-4}_2\CC^2$.
This gives the second direct summand of \eqref{U2onH12}, and we have finished a proof of
Proposition \ref{U2onH1}.
\proofend

\begin{remark}\label{rmk:zw}
{\em
From the above proof, elements of the representation spaces
$S^{n-2}_1\CC^2$ and $S^{n-4}_2\CC^2$ are polynomials in $u$.
As we have put $u=w/z$ where $(z,w)$ is the coordinates
on $\CC^2$, this is equivalent to saying that the representation spaces
are homogeneous polynomials of the two variables $z$ and $w$
(of degree $(n-2)$ and $(n-4)$ respectively.)
This will be useful later when identifying subgroups of $\U(2)$. 
}
\end{remark}

\subsection{Dimension of the moduli spaces}
Next we would like to compute,
by utilizing Proposition \ref{U2onH1}, dimension of the moduli space of ASD structures on the orbifold $\widehat{\ms O(-n)}$ obtained as small deformations of the LeBrun metric. 
For this, we need to 
compute dimension of orbits for the U(2)-action 
obtained in 
Proposition \ref{U2onH1}.

The case $n=3$ is very simple,
because the U(2)-action \eqref{U2onH12}
 is just a 2-dimensional representation $S^1_1\CC^2$.
If we restrict this to the subgroup SU(2), we get a natural representation of SU(2),
and any orbit  is
 diffeomorphic to a 3-sphere, except the origin.
Moreover U(2)-orbits and SU(2)-orbits
evidently coincide, and the stabilizer subgroup at any point
(except the origin) is isomorphic to U(1).


For investigating the case $n>3$, we
next compute dimension of U(2)-orbits in the space $S^m_l\CC^2$
for any $m\ge 2$ and $l\ge 0$.
We identify $S^m_l\CC^2$ with the space of 
homogeneous polynomials of $z$ and $w$ of degree $m$;
in particular a natural basis is provided by 
\begin{align}\label{monom}
z^m, z^{m-1}w, z^{m-2}w^2,\cdots, w^m.
\end{align}
Note that by Remark \ref{rmk:zw}, the variables $z,w$ are identical to the ones in 
the coordinates $(z,w)$ we used in Section \ref{ss:U2}.
We can classify all lower-dimensional orbits as follows:
%

\begin{proposition}\label{prop:orbit}
Suppose $m\ge 2$ and $l\ge 0$, and for each integer $j$ with
$0\le j\le m$, let $O_j\subset S^m_l\CC^2$ be the 
$\U(2)$-orbit going through the monomial $z^{m-j}w^j$.
Then we have the following:
(i) the orbit $O_j$ is  3-dimensional 
for any $j$,
(ii) the orbits $O_0,O_1,\cdots,O_m$ are all  $\U(2)$-orbits
in $S^m_l\CC^2$ which are not 4-dimensional,
except the origin,
(iii) the coincidence $O_j=O_{k}$ occurs iff 
$j=k$ or $j+k=m$ holds.

\end{proposition}

\proof
By thinking $(z,w)$ as homogeneous coordinates on 
$\CP^1$ we identify the space of homogeneous polynomials
of degree $m$ with $H^0(\CP^1,\ms O(m))$.
There is a natural $\U(2)$-action on this space, under which 
it is identified with $S^m\CC^2$
as a $\U(2)$-module.
Hence the $\U(2)$-module $S^m_l\CC^2$ is identified
with $\CC_l\otimes H^0(\ms O(m))$.
In order to classify lower-dimensional orbits, it is enough to 
classify all polynomials $P(z,w)\in \CC_l\otimes H^0(\ms O(m))$
whose stabilizer subgroup is of positive dimension.
We assert that this is the case exactly when 
the set of roots $Z_P:=\{(z,w)\in\CP^1\set P(z,w)=0\}$  
satisfies one of the following conditions:
(1)
$Z_P$ consists of a single point,
(2)
$Z_P$ consists of two points and moreover they are 
invariant under the involution 
$(z,w)\mapsto (\ol w, -\ol z)$ on $\CP^1$.

For this  suppose first that the polynomial
$P(z,w)$ satisfies the condition (1).
Then since the U(2)-action on $\CP^1$ is transitive,
we can suppose that $P(z,w) = az^m$ for some $a\in \CC^*$.
It is elementary to see that the
identity component of the stabilizer subgroup
for this monomial 
(viewed as an element of $ S^m_l\CC^2$) is a U(1)-subgroup
of $\U(2)$.
Therefore the stabilizer subgroup is 
of positive dimension.
Second suppose that 
$P(z,w)$ satisfies the condition (2).
Then under an identification
$\CP^1\simeq S^2$  the two roots of $P(z,w)=0$ (in $\CP^1$) 
form an anti-podal pair.
Recalling that 
the natural $\U(2)$-action on $\CP^1\simeq
S^2$ is isometric with respect to the standard metric on $S^2$,
this means that the natural U(2)-action
on the space of  anti-podal pairs
of points is also transitive.
Therefore we can suppose that 
$P(z,w)= z^{m-j} w^j$ for some
$0<j<m$.
Then again it is elementary to see that 
the stabilizer subgroup at $P(z,w)\in
S^m_l\CC^2$ is 1-dimensional.
Thus if $P(z,w)$ satisfies (1) or (2), then the stabilizer subgroup 
at this $P(z,w)\in S^m_l\CC^2$ is of 
1-dimensional.
This means the assertion (i)
of the proposition.

Conversely for (ii) suppose that the roots of $P(z,w)\in S^m_l\CC^2$ do not satisfy (1) nor (2).
If there are more than two roots, then elements of 
$\U(2)$ preserving the set of the roots constitute a finite 
subgroup at most. This implies that 
the stabilizer subgroup at the point $P(z,w)\in S^m_l\CC^2$ is 
also a finite subgroup.
If there are exactly two roots but  the roots are not an
anti-podal pair, then because of the isometricity of the 
$\U(2)$-action on $S^2$, elements of $\U(2)$ which preserve
the two roots constitute a finite subgroup at most,
since such an element has to preserve four points.
Therefore we again obtain that the stabilizer subgroup
at $P(z,w)\in S^m_l\CC^2$ is a finite subgroup.
Therefore the stabilizer subgroup is
zero-dimensional if $P(z,w)$ does not satisfy (1) nor (2).
Hence from the transitivity of 
the natural U(2)-action on $S^2$,
we obtain that 
if the U(2)-orbit through $P(z,w)\in S^m_l\CC^2$ is not four-dimensional,
then $P(z,w)\in O_j$ for 
some $0\le j\le m$.
This proves the assertion (ii).

For the final assertion (iii), 
$O_j=O_{m-j}$ is clear  because there actually exists
an element of U(2) which interchanges $z^{m-j}w^j$
and $z^jw^{m-j}$ as elements of $S^m_l\CC^2$.
Moreover, if $O_j=O_k$, the set of multiplicities of the 
two roots must equal. 
This implies
$k\in \{j,m-j\}$.
\proofend

\vsp
The stabilizer subgroup at the monomials 
in the space \eqref{U2onH12} is concretely given as follows: 

\begin{lemma}\label{lemma:stab}
For a pair $(m_1,m_2)$ of integers define a subgroup
$G(m_1,m_2)\subset T^2\subset U(2)$ by 
\begin{align}
G(m_1,m_2):= \left\{
\begin{pmatrix}
e^{i\aaa} & 0 \\
0 & e^{i\bbb}
\end{pmatrix}
; \aaa,\bbb\in\RR, \,m_1\aaa+ m_2\bbb=0
\right\}.
\end{align}
Then the identity component of the stabilizer subgroup
at the monomial $z^{m-k}w^k\in S^m_l\CC^2$ is 
$G(m+l-k,l+k)$. 
In particular the identity component coincides
for the two monomials 
\begin{align}
z^{n-2-k}w^k\in S^{n-2}_1\CC^2\,
{\text{ and }} \,z^{n-3-k}w^{k-1}\in S^{n-4}_2\CC^2
\end{align}
in the representation \eqref{U2onH12}.
\end{lemma}

\proof
This is elementary and we omit a proof.
The last coincidence is a direct consequence of the
concrete form of the stabilizer subgroup.
\proofend

\vsp
From Lemma \ref{lemma:stab} we put all monomials in $S^{n-2}_1\CC^2$ and 
$S^{n-4}_2\CC^2$ as in the following table:

\vsp
\begin{tabular}{c| ccccccc}
${\small {S^{n-2}_1\CC^2}}$ & $z^{n-2}$ & $z^{n-3}w$ & $z^{n-4}w^2$ & 
$\cdots$ & $z^2 w^{n-4}$ & $zw^{n-3}$ & $w^{n-2}$ \\
${\small S^{n-4}_2\CC^2}$ & & $z^{n-4}$ & $z^{n-5}w$ &
 $\cdots$ & 
$zw^{n-3}$ & $w^{n-4}$\\
stabilizer & {\small $G(n-1,1)$} & 
{\small $G(n-2,2)$} & {\small $G(n-3,3)$}
& $\cdots$ & {\small $G(3,n-3)$} & 
 {\small $G(2,n-2)$} &  {\small $G(1,n-1)$}
\end{tabular}

\vsp\noindent
This reads, for example, that the two monomials
$z^{n-3}w$ and $z^{n-4}$ have the same stabilizer group
$G(n-2,2)$ as the identity component.

With these preliminary results,
we investigate automorphisms and the moduli space of the ASD structures
on $\widehat{\ms O(-n)}$ which appear from the versal family 
of the twistor space $Z_{\LB}$.
Let $p:\ms Z \to B$ be the versal family
for $Z_{\LB}$,
where $p^{-1}(0)=Z_{\LB}$.
As before, the parameter space
$B$ can be U(2)-equivariantly
identified with an invariant neighborhood of the origin in
the vector space $S^{n-2}_1\CC^2\oplus S^{n-4}_2\CC^2$.
For each subgroup $G\subset\U(2)$
let $B^G$ the subspace of $G$-invariant elements,
which is an intersection of $B$
with the linear subspace of $S^{n-2}_1\CC^2\oplus S^{n-4}_2\CC^2$ consisting of $G$-invariant elements.
By restricting $p$
over  $B^G$,
we obtain a versal family of $G$-equivariant deformations
for $Z_{\LB}$.
In particular $Z_{\LB}$ admits
a non-trivial  $G$-equivariant deformation
iff $B^G\neq 0$.
These considerations readily mean the following 

\begin{proposition}\label{prop:edf}
Suppose $n\ge 3$ and let $G$ be a closed connected subgroup of $\U(2)$ which satisfies $\dim G\ge 1$.
Assume that
the LeBrun's ASD structure on $\widehat{\ms O(-n)}$ admits
a non-trivial $G$-equivariant deformation.
Then (i) if $n\neq 4$, 
we have $G=G(k,n-k)$ for some $k$ satisfying 
$1\le k<n$.
So $G$ is isomorphic to $\U(1)$.
(ii) If $n=4$ and $G\neq \SU(2)$,
the same conclusion holds.
Moreover there exists an $\SU(2)$-equivariant deformation.
\end{proposition}

\proof
The assertion except the final one in (ii)
follows immediately from
Propositions \ref{U2onH1}, \ref{prop:orbit}
and the above table.
When $n=4$, the second direct summand in \eqref{U2onH12} becomes the 1-dimensional space $\CC_2$,
and the identity component of the stabilizer
subgroup at any point on $\CC_2$ is clearly the subgroup SU(2),
except the origin.
Therefore if we restrict the versal family 
of ASD structures on $\widehat{\ms O(-4)}$ to the
real 2-dimensional subspace
$\{0\}\oplus \CC_2\subset S^2_1\CC^2\oplus\CC_2\simeq H^1
(\Theta_{Z_{\LB}})^{\sigma}$, we obtain the required SU(2)-equivariant deformation.
\proofend

\vsp
For the moduli space of the invariant ASD structures
in Proposition \ref{prop:edf}, if $t\in B^G$,
 the fiber $p^{-1}(t)$ has a (holomorphic)
$G$-action of course.
However, there can exist a subgroup of U(2)
which acts non-trivially on $B^G$,
and it gives
an identification between different fibers of $p$
over $B^G$.
Thus
the subspace $B^G$ itself cannot be considered as a moduli space
of $G$-invariant ASD structures in general,
and instead
the actual moduli space is considered to be the quotient space of 
the subspace $B^G$ under the action of 
the subgroup of  U(2) consisting of elements
which preserve $B^G$.

For example, if $G=G(n-2,2)$,
the subspace $B^G$
is $B\cap \langle z^{n-3}w,z^{n-4}\rangle_{\CC}$,
which is 4-dimensional over $\RR$.
It is easy to see that the subgroup of U(2) consisting of elements which 
preserve this subspace is the maximal torus $T^2$
consisting of diagonal matrices.
Further orbits of the $T^2$-action on 
$\langle z^{n-3}w,z^{n-4}\rangle_{\CC}$ is 1-dimensional,
except the origin.
Consequently we obtain that the moduli space of 
these U(1)-invariant ASD structures
on $\widehat{\ms O(-n)}$ is $3$-dimensional.
By the same argument, we get the following

\begin{proposition}\label{prop:sub1}
Let $n\ge 3$ and $k$ satisfy $1\le k<n$, and 
consider  ASD structures on $\widehat{\ms O(-n)}$
obtained as the $G(k,n-k)$-equivariant small deformation
as in Proposition \ref{prop:edf}.
Then the moduli space of $\U(1)$-invariant ASD structures
on $\widehat{\ms O(-n)}$ obtained by this equivariant
deformation is $1$-dimensional if $k\in \{1,n-1\}$, and $3$-dimensional if $k\not\in \{1,n-1\}$
For the $\SU(2)$-equivariant deformation in the case $n=4$,
the moduli space is 1-dimensional.
\end{proposition}

It is already immediate to give a proof of Theorem \ref{thm:smry1} in the introduction.

\vsp
\noindent
{\em Proof of Theorem \ref{thm:smry1}.}
This is an immediate consequence of 
Propositions \ref{U2onH1}, \ref{prop:orbit}, \ref{prop:sub1}
and Lemma \ref{lemma:stab}.
More concretely for the U(2)-invariant subset $B_i$ in
$B$ in the theorem, it is enough to take the union of all U(2)-orbits
which go through 
\begin{itemize}
\item the real 2-dimensional subspace $\langle z^{n-2} \rangle_{\CC}\backslash\{0\}$ for the 
case $i=1$,
\item the real $4$-dimensional subspace
$\langle z^{n-i-1}w^{i-1}, z^{n-i-2}w^{i-2}\rangle_{\CC}
\backslash\{0\}$ for the case $1<i\le [n/2]$,
\item the real $2$-dimensional subspace 
$(0\oplus \CC_2)\backslash \{0\}$ in the case 
of $(n,i)=(4,0)$.\proofend
\end{itemize} 

\vspace{1mm}
For the dimension of the above $\U(2)$-invariant subsets,
we readily have
$\dim B_0=2,$ $\dim B_1 = 4$ and $\dim B_i= 6$
if $i\not\in\{0,1\}$.

\section{Deformations preserving K\"ahlerian property}
\label{s:sfk}
The investigation in the last section concerns
versal and  equivariant deformations of the LeBrun's ASD structure 
on $\widehat{\ms O(-n)}$ as an ASD orbifold.
Since the LeBrun metric is K\"ahler,
from differential geometric point of view, it would be
desirable to obtain deformations of the metric preserving 
not only anti-self-duality but also K\"ahlerian property.
In this section again by using twistor spaces, we find,
for any $n\ge 3$, a deformation of the LeBrun metric
on $\ms O(-n)$
which keeps anti-self-duality as well as K\"ahlerity.
The deformation is realized as one of the
U(1)-equivariant deformation
we found in the last section.
Meanwhile we also show that 
the corresponding twistor spaces
are Moishezon for these deformations.
We also show that for other U(1)-equivariant
deformations, the deformed twistor spaces
are not Moishezon.

The key tool for finding such a deformation is of course Pontecorvo's theorem
\cite[Theorem 2.1]{Pont92}, which means that an anti-self-dual conformal structure on a 4-manifold $M$ carries a K\"ahler representative
for a complex structure if and only if the twistor space
possesses a divisor $D$ which is mapped diffeomorphically
to $M$ by the twistor fibration, and which satisfies $D+\ol D\in |K^{-1/2}|$; then the conformal class has
a K\"ahler representative with respect to the complex structure of $D$, and then derive information about
existence of a reducible member of $|K^{-1/2}|$.

As we already mentioned,
for the twistor space $Z_{\LB}$ of the 
LeBrun orbifold $\widehat{\ms O(-n)}$,
the divisor $D$  in
Section \ref{ss:recall} gives a reason
for
the LeBrun metric to be K\"ahler with respect to the standard complex structure on $\ms O(-n)$.
One would naturally think that 
we should investigate deformations of 
the pair $(Z_{\LB},D+\ol D)$,
which might actually work.
However, the divisor $D$ itself is not a Cartier divisor
on $Z_{\LB}$,
and deformation theory of such a pair
might be subtle.
Therefore here we take  
a real irreducible divisor $S\in |K^{-1/2}|$
and consider deformations of the pair
$(Z_{\LB},S)$.

\subsection{Deformation of the pair $(Z_{\LB},S)$.}
\label{ss:dfmpair}
In general if $X$ is a complex variety 
and $Y$ is a reduced Cartier divisor on $X$, 
the subsheaf $\Theta_{X,Y}$ of the tangent sheaf
$\Theta_X$ is naturally defined as
\begin{align}
\Theta_{X,Y} = \{v\in\Theta_X\set v(g)/g\in \ms O_X\},
\end{align}
where $g\in \ms O_X$ is a local equation of $Y$.
(When $X$ and $Y$ are non-singular, this is exactly
the sheaf of vector fields on $X$ which are tangent to $Y$.)
This sheaf plays the same role for deformations
of the pair $(X,Y)$ as the sheaf $\Theta_X$ plays for deformations
of $X$ itself.
Namely, if $X$ is compact,
the cohomology group $H^1(\Theta_{X,Y})$ 
is the Zariski tangent space of the Kuranishi family of {\em locally trivial} deformations of the pair $(X,Y)$,
and $H^2(\Theta_{X,Y})$ is the obstruction space.
In particular, if $H^2(\Theta_{X,Y})=0$,
the parameter space of the Kuranishi family 
is naturally identified with an open neighborhood
of the origin in $H^1(\Theta_{X,Y})$.

As a real irreducible divisor $S\in |K^{-1/2}|$ we first take any 
real non-singular $(1,1)$-curve
$\ms C$ on $Q=\CP^1\times\CP^1$
which is different from the diagonal $\Delta$,
and put $S:=\mu(\pi^{-1}(\ms C))$.
(Recall that $\pi:X\to Q$ is a projection and $\mu:X\to Z_{\LB}$ is the blowdown of the divisor $E\sqcup \ol E$.)
From the formula \eqref{cbl},
this  actually belongs to $|K^{-1/2}|$ on $Z_{\LB}$.
Then the maximal subgroup of the automorphism group U(2)
of $Z_{\LB}$ which preserves $S$ is isomorphic to a torus
$T^2$, and under this action 
$S$ has a structure of a toric surface.
We note that the complex structure of $S$ is independent
of the choice of the $(1,1)$-curve $\ms C$.
Obviously the above subgroup $T^2\subset \U(2)$
preserves not only $S$ but also the two divisors $D$ and $\ol D$ too.
It is easy to see that the surface $S$ satisfies the following properties:
\begin{itemize}
\item the intersection $S\cap L_{\infty}$ consists of
two  points, which are mutually conjugate, 
\item $S$ has $A_{n-1}$-singularities at these two points,
and is non-singular except these points.
\end{itemize}

For investigating deformation of the pair
$(Z_{\LB},S)$, we first show the following

\begin{proposition}\label{prop:van1}
We have the following vanishing:
\begin{align} 
H^2(\Theta_S) = 
H^2(\Theta_{Z_{\LB}}(-S)) = 0.
\end{align}
\end{proposition}

\proof
In this proof we
again write $Z$ for $Z_{\LB}$.
Our proof for  $H^2(\Theta_S)=0$ is quite analogous to 
$H^2(\Theta_Z)=0$ in the proof of Proposition \ref{prop:unobs}.
Put $S\cap L_{\infty}=\{p,\ol p\} $.
Since $S$ is biholomorphic 
to $\pi^{-1}(\ms C)\,(\subset X)$,
we have two exact sequences
\begin{align}\label{ses:6}
0\lras \Theta_{S/\ms C} \lras \Theta_S
\lras \ms G\lras 0 \,\,{\text{ and }}\,\,
0 \lras\ms G\lras \pi^*\Theta_{\ms C}
\lras \CC_p\oplus \CC_{\ol p}
\lras 0,
\end{align}
where $\ms G$ is the image sheaf of the canonical homomorphism
$\Theta_S\to\pi^*\Theta_{\ms C}$.
Since $\ms C\simeq\CP^1$, 
the induced map $H^0(\pi^*\Theta_{\ms C})\to H^0(\CC_p\oplus
\CC_{\ol p})$ is easily seen to be surjective.
Further we have  $H^i(\pi^*\Theta_{\ms C})\simeq H^i(\Theta_{\ms C})=0 $
for  any $i\ge 0$.
Hence the second sequence of \eqref{ses:6} implies
$H^i(\ms G)=0$ for $i\ge 1$.
Hence the first one in \eqref{ses:6} means $H^2(\Theta_{S/\ms C})\simeq H^2(\Theta_S)$.
Further, from an obvious vector field
which vanishes on  $C\sqcup\ol C$, we have
$\Theta_{S/\ms C}\simeq\ms O_S(C+ \ol C)$,
and we readily have $H^2(\ms O_S(C + \ol C))=0$.
Therefore $H^2(\Theta_S)=0$ follows.

In the sequel we put $F:=K^{-1/2}$ for simplicity,
and show $H^2(\Theta_Z\otimes F^{-1})=0$.
By taking tensor product with $F^{-1}$ 
to the exact sequence \eqref{ses1}, 
we obtain an exact sequence
\begin{align}\label{ses2}
0 \lras \Theta_{Z,C+\ol C}\otimes F^{-1} \lras
\Theta_Z \otimes F^{-1}\lras 
(N_{C/Z}\oplus N_{\ol C/Z})\otimes F^{-1}|_C
\lras 0.
\end{align}
Further since $F|_C\simeq K_S^{-1}|_C\simeq
\ms O_C(2-n)$ and $N_{C/Z}\simeq
\ms O_C(1-n)^{\oplus 2}$, 
the last non-trivial sheaf in \eqref{ses2} is 
isomorphic to $\ms O_C(-1)^{\oplus 2}
\oplus \ms O_{\ol C}(-1)^{\oplus 2}$.
Hence we have 
\begin{align}\label{iso0}
H^2(\Theta_Z\otimes F^{-1})\simeq
H^2(\Theta_{Z,C+\ol C}\otimes F^{-1}).
\end{align}
For computing the RHS, 
from the isomorphism
$
\Theta_{X,E+\ol E}\simeq\mu^*\Theta_{Z,C+\ol C},
$
recalling $\mu^*F\simeq\pi^*\ms O_Q(1,1)\otimes \ms O_X(E+\ol E)$, we have $\mu^*(\Theta_{Z,C+\ol C}\otimes F^{-1})\simeq\Theta_{X,E+\ol E}\otimes
\pi^*\ms O_Q(-1,-1)\otimes \ms O_X(-E-\ol E)$.
For simplicity we write $\ms L$ for the sheaf on RHS.
From the last isomorphism we have
\begin{align}\label{iso1}
H^2(\Theta_{Z,C+\ol C}\otimes F^{-1})
\simeq
H^2(X, \ms L).
\end{align}
For the RHS of this, from the inclusion
$0\to \Theta_{X,E+\ol E}\to \Theta_X$
we have the exact sequence
\begin{multline}
\label{ses3}
0 \lras
\ms L 
\lras
\Theta_{X}\otimes
\pi^*\ms O_Q(-1,-1)\otimes \ms O_X(-E-\ol E)\\
\lras
(N_{E/X}\otimes \ms O_X(-E)
\oplus
N_{\ol E/X}\otimes \ms O_X(-\ol E))
\otimes
\pi^*\ms O_Q(-1,-1)
\lras 0.
\end{multline}
The last non-trivial term of \eqref{ses3} is
clearly isomorphic to $\pi^*\ms O_Q(-1,-1)
|_{E\sqcup\ol E}$,
whose all cohomologies vanish.
Therefore, writing the middle sheaf
as $\ms L'$, we get
\begin{align}\label{iso2}
H^2(\ms L) \simeq H^2(\ms L').
\end{align}
For the RHS of this, by taking  tensor
product with
$\pi^*\ms O_Q(-1,-1)\otimes \ms O_X(-E-\ol E)$  to the exact sequence \eqref{ses:v},
we obtain
\begin{align}\label{ses:4}
0\lras \pi^*\ms O_Q(-1,-1)
\lras \ms L' \lras
\ms F\otimes \pi^*\ms O_Q(-1,-1)\otimes \ms O_X(-E-\ol E) \lras 0.
\end{align}
Writing $\ms F'$ for the last non-trivial sheaf of this sequence, by taking a tensor product with
$\pi^*\ms O_Q(-1,-1)\otimes \ms O_X(-E-\ol E)$  to the exact sequence \eqref{ses0}, we get
\begin{align}
0\lras \ms F'
\lras 
\pi^*(\ms O(1,-1)\oplus\ms O(-1,1)) 
\lras
\ms O_{\Delta} \lras 0.
\end{align}
From this we get $H^2(\ms F')=0$.
Therefore from \eqref{ses:4} we obtain 
$H^2(\ms L')=0$.
Hence by \eqref{iso2}, \eqref{iso1} and \eqref{iso0} we obtain 
$H^2(\Theta_Z\otimes F^{-1})=0$.
\proofend

\vsp
For investigating deformations of the pair $(Z_{\LB},S)$
we also need the following
\begin{proposition}\label{prop:van2}
Let $S$ be a real irreducible member
of $|K^{-1/2}|$ as above. Then
we have the following exact sequence
\begin{align}\label{ses:5}
0\lras
\Theta_{Z_{\LB}}(-S)
\lras \Theta_{Z_{\LB},\,S} \lras 
\Theta_S
\lras 0.
\end{align}
Hence by Proposition \ref{prop:van1}, we have
\begin{align}  H^2(\Theta_{Z_{\LB},\,S})=0.
\end{align}
\end{proposition}

\proof
We again write $Z$ for $Z_{\LB}$.
Since $\Sing S\subset\Sing Z$,
\eqref{ses:5} is obvious outside the two singular points of $S$.
Also, in a neighborhood of the singular points, 
defining equation of $Z$ and $S$ in the ambient space
$\mathbb P(\ms E_n)$ can be taken as
$xy-(u-v)^n=0$ and $xy-(u-v)^n=u=0$ respectively, 
and by using these and \eqref{tansh} it is easy to obtain 
concrete form of sections of the sheaves
$\Theta_{Z,S}$ and $\Theta_S$
in a neighborhood of the singular point.
From this the exact sequence \eqref{ses:5}
is known to be available on the two singular points too.
%
\proofend

\vsp
From Propositions \ref{prop:van1} and \ref{prop:van2} 
we readily obtain the following co-stability:

\begin{proposition}
\label{prop:cost}The irreducible Cartier divisor $S$ is co-stable
in $Z_{\LB}$ 
with respect to locally trivial deformations of $S$.
Namely for any such deformation there exists a locally 
trivial deformation of the pair $(Z_{\rm LB}, S)$ which gives
the prescribed deformation of $S$ by restriction.
\end{proposition}

\proof
We again write $Z$ for $Z_{\LB}$.
Let $p:\mathscr Z\to B$ and $\ms S\subset \ms Z$ be the Kuranishi family of locally trivial deformations of the pair $(Z,S)$, 
where $p^{-1}(0)=Z$ and $p^{-1}(0)\cap \ms S=S$.
By Proposition \ref{prop:van2}, $B$ can be identified 
with an open neighborhood of $0$ in $H^1(\Theta_{Z,S})$. 
Let $\ms S'\to B'$ be the Kuranishi family of 
locally trivial deformation of $S$.
As $H^2(\Theta_S)=0$ by Proposition \ref{prop:van1},
the parameter space $B'$ may be  identified with
an open neighborhood of the origin in $H^1(\Theta_S)$.
By versality of the Kuranishi family,
the family $\ms S\to B$ induces a holomorphic map
$f:B\to B'$ with $f(0)=0$, and the differential $df$ at $0$ is
identified with the natural linear map $H^1(\Theta_{Z,S})\to H^1(\Theta_S)$. 
The last map is locally submersion by Proposition \ref{prop:van2}.
Therefore $f$ is locally surjective at $0$.
By the property $f^*\ms S'\simeq \ms S$ over $B$,
this means the required co-stability.  
\proofend

\subsection{Concrete deformations of the surface $S$, and deformations of the pair}
\label{ss:dfmS}
Next we concretely construct locally trivial deformations of the 
singular toric surface
$S$ which preserve  U(1)-action, for some explicit subgroups
$\U(1)$ in $\U(2)$.
Applying Proposition \ref{prop:cost}
to any one of these deformations,
we will obtain non-trivial deformations of the LeBrun metric.
It will turn out that some of these deformations
preserve K\"ahlerian property.

Fix any $n\ge 3$ as before. We first  realize our singular toric  surface $S$
in $Z_{\LB}$
as an explicit birational transform of 
 the product surface $\CP^1\times\CP^1$.
Writing $0:=(1:0)\in\CP^1$ and $\infty:=(0:1)\in\CP^1$,
We take four points and four curves on $\CP^1\times\CP^1$ as 
\begin{gather*}
q_1=(0,0),\,q_2=(\infty,0),\,q_3=(0,\infty),\,
q_4=(\infty,\infty),\\
C_1=\CP^1\times 0, \,C_2 = \infty \times\CP^1,\,
C_3=\CP^1\times \infty,\,C_4=0\times\CP^1.
\end{gather*}
We regard
$\CP^1\times\CP^1$ as a toric surface
by considering the  $T^2$-action which preserves the curve $C_1+C_2+C_3+C_4$.

For any integer $k$ satisfying $0< k < n$,
we assign a weight $k$ on the two points $q_1$ and $q_3$,
and 
a weight $(n-k)$ on the other points $q_2$ and $q_4$.
Under this setting let $\tilde S\to \CP^1\times \CP^1$
be the surface obtained by blowing-up the  point
$q_i$ for $m_i$ times for any $1\le i\le 4$, where $m_i$ is the above weight at $q_i$.
Here, if $m_i\ge 2$, the blowup is always done at a $T^2$-fixed point on
the strict transforms of
the curve $C_1$ or $C_3$.
$\tilde S$ is also a toric surface.
The inverse image of the curve
$C_1+C_2+C_3+C_4$ is a
 $T^2$-invariant anticanonical curve on $\tilde S$,
and it
consists of $4+2n$ components.
The self-intersection numbers of the components are
given by, up to cyclic permutations,
\begin{align}\label{ac1}
-n,-1,
\overbrace{-2,\cdots,-2}^{n-1},-1,-n,-1,\overbrace{-2,\cdots,-2}^{n-1},-1.
\end{align}
In particular, these are independent of $k$,
and therefore so is the structure of the toric surface $\tilde S$.
(Dependence on $k$ will appear later.)
Let $C$ and $\ol C$ be the two $(-n)$-curves among \eqref{ac1}.
These are strict transforms of the curves $C_1$
and $C_3$.
Then in $\tilde S$ we can contract the two chains of the $(-2)$-curves to
obtain a toric surface with two $A_{n-1}$-singularities.
By looking structure as a toric surface, it is easy to see that the last surface is biholomorphic to the surface $S$
in $Z_{\LB}$ we have given in the previous subsection.
The contraction map $\tilde S\to S$ is nothing but the minimal resolution 
of the singularities of $S$,
and $C$ and $\ol C$ are exactly the curves $\mu(E)$ and $\mu(\ol E)$ under the identification.
It is also easy to verify that if we introduce a real structure 
on the initial surface $\CP^1\times\CP^1$ by the product which is (anti-podal)$\times$
(complex conjugation), then it naturally lifts 
to be a real structure on the surface $\tilde S$ as well as that on the contracted surface,
and the last real structure is exactly the one
on the real divisor $S$ in $Z_{\LB}$.

Now we shall give U(1)-equivariant, locally trivial deformations
of the surface $S$ preserving the real structure,
by using the above realization of $S$.
The deformations we  construct 
are uniquely and explicitly determined
from the value $k$ above.
In order to construct  locally trivial deformation of $S$,
it is enough to give a deformation (in the usual sense) of the minimal resolution 
$\tilde S$ which preserves the two chains of $(-2)$-curves.
For fixed integer $k$ with $0<k<n$ as above, we think the surface $\tilde S$ as obtained 
by blowing up $\CP^1\times\CP^1$ 
in the way indicated by
the weights $k$ and $n-k$ as above.
Then by moving the weighted blowup points $q_1$ and $q_3$
along the curves $C_4$ and $C_2$ freely respectively,
we obtain a 2-dimensional family of smooth rational surfaces which can naturally be regarded as
deformation of the
surface $\tilde S$.
If $m_i$ is the weight at the point $q_i$ 
as above, even after the deformation,
the iterated blowups at $q_i$ yield $(m_i-1)$ number
of $(-2)$-curves as exceptional curves.
Also, the strict transforms of the two curves
$C_2$ and $C_4$ are still $(-2)$-curves in the deformed new surface.
The union of all these $(-2)$-curves still form two chains of $(-2)$-curves, and  each chain yet consists of $(n-1)$ components.
Hence by contracting these two chains
simultaneously, 
we obtain a family of 
rational surfaces which have two $A_{n-1}$-singularities.

In this way, for each $1\le k<n$ we have obtained a locally trivial deformation
of the toric surface $S$.
From the construction the parameter space of this
deformation is naturally identified with the 
product $C_2\times C_4$.
The real structure on $\CP^1\times\CP^1$ attached above
naturally acts on this product, and 
by restricting the deformation to the real locus,
we obtain a deformation of $S$ preserving the real structure.
The parameter space of this family is clearly real 2-dimensional.
These deformations actually deform the complex structure
of $S$ since the deformed new surface is not a toric surface anymore.

Next we show that all these deformations
(determined by $k$) of the surface $S$ are equivariant
with respect to a U(1)-subgroup of $T^2$,
and the subgroup depends on the value $k$.
For this we consider the U(1)-action
on $\CP^1\times\CP^1$ which fixes points
on $C_2\cup C_4$.
This U(1)-action
clearly fixes  the four points $q_1,\cdots,q_4$, even after moving $q_1$ and $q_3$.
Moreover, it is immediate to see that 
 the U(1)-action lifts on the blowup
even after moving,
and that on the final surface
the chains of $(-2)$-curves are  
invariant under the induced U(1)-action.
Therefore the U(1)-action descends 
on the contraction of the two chains.
Thus our deformation of the surface $S$
 is  U(1)-equivariant.

This U(1) can be naturally regarded as a subgroup
of the torus $T^2$, where the last $T^2$ is thought
as an automorphism group of the toric surface $\tilde S$
preserving the real structure.
Then this subgroup has to depend on the number $k$,
since the component fixed by the subgroup depends on $k$.
(More concretely, there are exactly $k$ components
between the curve $C$ and the fixed component.)
We write $G(k)$ for this U(1)-subgroup of $T^2$.
While $G(k)$ is a subgroup of $T^2$,
it is naturally regarded as a subgroup of the automorphism
group U(2) of the LeBrun twistor space $Z_{\LB}$,
since the torus $T^2$ is originally the
maximal subgroup of 
U(2) which preserves the divisor $S$.

On the other hand in Proposition \ref{prop:edf}
for each $1\le k<n$ we have
obtained the U(1)-subgroup $G(k,n-k)$ for which
$Z_{\LB}$ admits a non-trivial equivariant deformation.
These subgroups coincide:

\begin{proposition}\label{prop:identify}
For any $1\le k<n$, we have $G(k,n-k)=G(k)$ in $\U(2)$.
\end{proposition}

\proof
For the divisors $S$ and $D+\ol D$ of $|K^{-1/2}|$
on $Z_{\LB}$, 
the intersection 
$S\cap (D\cup\ol D)$ consists of a cycle of six smooth rational curves, two of which are $C$ and $\ol C$, while the remaining four components
are the inverse image of the two intersection points
$\Delta\cap \ms C\subset Q$
 under the projection $\pi$.
This cycle is naturally divided into halves by 
the twistor line $L_{\infty}$.
As $S$ and $D\cup\ol D$ are invariant under the $T^2$-action, this cycle is also $T^2$-invariant.
Moreover the $T^2$-action on the cycle is effective.
Therefore we can identify any U(1)-subgroup of $T^2$ 
from the action on each component of the cycle.

For the subgroup $G(k,n-k)$, we can readily obtain
these actions in a concrete form,
by recalling that the torus $T^2$ in which $G(k,n-k)$ is included
is exactly the maximal torus of $\U(2)$ which consists
of diagonal matrices, and that the U(2)-action on $\ms O(-n)$ can be obtained naturally via 
$\ZZ_n$-quotient and the minimal resolution.
(Recall also that the divisor $D\backslash L_{\infty}$
in $Z_{\LB}\backslash L_{\infty}$ is U(2)-equivariantly
identified with the open subset $\ms O(-n)$
by the twistor fibration map.)
On the other hand, the $T^2$-action on $S$ was also explicitly 
given through the above construction, 
and therefore we can easily recognize the action
of the subgroup $G(k)$ on the cycle in concrete forms.
The coincidence $G(k,n-k)=G(k)$ follows from these explicit computations.
We omit the detail.
\proofend

\vsp
From the proposition, by making use of the co-stability obtained in Proposition \ref{prop:cost} we now have the following

\begin{proposition}\label{prop:dfmpair}
For each integer $k$ with $1\le k<n$, let $G(k,n-k)\subset\U(2)$ be the $\U(1)$-subgroup  given
as in Proposition \ref{prop:sub1},
and $S$  the irreducible member of $|K^{-1/2}|$ as 
taken in Section \ref{ss:dfmpair}.
Then the pair $(Z_{\LB},S)$ admits a
$G(k,n-k)$-equivariant, locally trivial deformation.
Moreover, the twistor space $Z_{\LB}$ itself actually deforms
in this deformation.
\end{proposition}

\proof
As constructed above, the surface $S$ admits a
$G(k)$-equivariant deformation for which the complex 
structure actually deforms.
By Proposition \ref{prop:identify}, 
this deformation is also $G(k,n-k)$-equivariant.
Applying
 Proposition \ref{prop:cost} to this deformation of $S$, 
there exists a locally trivial 
deformation of the pair $(Z,S)$ 
which induces the last deformation of $S$ by restriction,
where $Z=Z_{\LB}$ as before.
This deformation of the pair can be taken $G(k,n-k)$-equivariantly, since the exact sequence
\eqref{ses:5} is $T^2$-equivariant,
so that the induced map $H^1(\Theta_{Z,S})\to
H^1(\Theta_S)$ is also $T^2$-equivariant.
Thus we obtain the existence of the $G(k,n-k)$-equivariant deformation of the pair $(Z,S)$.
The complex structure of $Z$ actually varies in
this deformation, since
the complex structure of the divisor $S$ actually deforms,
while when we move $S$ inside $Z$, 
$S$ remains to be a toric surface, so that the complex structure does
not vary.
\proofend

\vsp
From the proposition, the $G(k,n-k)$-equivariant deformation
found in Proposition \ref{prop:sub1}
of the LeBrun orbifold
 can be
realized by a deformation of the twistor space
for which the divisor $S$ survives.
However, since we are taking an
irreducible $S\in |K^{-1/2}|$ 
and not taking the reducible divisor $D+\ol D$,  
we do not know at this stage if the metric
is accordingly deformed in a way that 
the K\"ahlerian property with respect to some complex structure on $\ms O(-n)$ is preserved.
In the next subsection we answer this affirmatively for $k\in \{1,n-1\}$.

\subsection{Deformation preserving K\"ahlerian property}
For that purpose we first investigate
 pluri-anticanonical systems of the singular rational surfaces
obtained by
the $G(k)$-equivariant deformation in Section \ref{ss:dfmS}.
As in the previous subsection 
let $n\ge 3$ and $S$ be the real irreducible member of $|K^{-1/2}|$ on $Z_{\LB}$ 
for the LeBrun metric on $\widehat{\ms O(-n)}$.
We write $S_t$ for the singular rational
surface obtained 
by the $G(k)$-equivariant deformation
of $S$ constructed in the last subsection.
Since $S_t$ also has only $A_{n-1}$-singularities,
the canonical divisor $K$ of $S$ is a Cartier divisor, so the system $|mK^{-1}|$ and 
the anti-Kodaira dimension
(i.e.\,the Kodaira dimension of $K^{-1}$) makes sense.
We denote the latter by $\kappa^{-1}(S_t)$
as usual.

\begin{proposition}
\label{prop:dfmS}
The rational surface
$S_t$ satisfies the following properties:
(i) if $k\in \{1,n-1\}$,
 we have $\kappa^{-1}(S_t)=2$,
(ii) if $n=4$ and $k=2$, we have $\kappa^{-1}(S_t)=1$,
(iii) if $n>4$ and $k\not\in\{1,n-1\}$,
we have $|mK^{-1}|=\emptyset$ for any $m>0$.
\end{proposition}

\proof
For (i), from the construction of the
$G(k)$-equivariant deformation,
the surface $S_t$ is obtained from 
a non-singular surface 
by contracting two chains of $(-2)$-curves.
We denote the last non-singular surface
by $\tilde S_t$.
(So the contraction $\tilde S_t\to S_t$
is the minimal resolution.)
If $k\in\{1,n-1\}$, the  surface $\tilde S_t$ is exactly  the divisor in $|K^{-1/2}|$ on the twistor space over 
$n\CP^2$ that we have investigated in \cite{Hon07-2,
HonDSn1}.
In particular, the system $|(n-2)K^{-1}|$ on 
$\tilde S_t$ induces a surjective degree-two morphism 
 $\tilde S_t\to\CP^2$.
Hence, since $S_t$ has only $A_{n-1}$-singularities, the degree-two morphism $\tilde S_t\to \CP^2$ factors as 
$\tilde S_t\to S_t\to \CP^2$, where 
$\tilde S_t\to S_t$ is the contraction 
and $S_t\to\CP^2$ is the map associated to
$|(n-2)K^{-1} |$ on $S_t$.
This implies $\kappa^{-1}(S_t)=2$.

For (ii) let $\tilde S_t\to S_t$ have the same meaning as above.
Then if $n=4$ and $k=2$, 
the surface $\tilde S_t$ is 
the same as the divisor in $|K^{-1/2}|$ on the twistor spaces on $4\CP^2$ of algebraic dimension  two which was investigated in \cite{Hon01}.
In particular $|K^{-1}|$ on $\tilde S_t$ is base point free and induces an elliptic fibration $f:\tilde S_t\to\CP^1$.
Hence we have $f^*\ms O(1)\simeq K^{-1}$,
which means 
$\kappa^{-1}(\tilde S_t)=1$.

For (iii) we first consider the surface $\tilde S_t$ 
in the case $n=5$ and $k\in \{2,3\}$.
Obviously this surface is obtained from
the elliptic surface in the last case
of $(n,k)=(4,2)$ by blowing-up two points.
Further the two points belong to
mutually different fibers of the elliptic fibration, and moreover the two points are
smooth point of the fibers. From these we readily deduce that
 $h^0(mK^{-1})=0$ on $\tilde S_t$ for any $m>0$. 
Therefore since $h^0(mK^{-1})$ for fixed $m$ cannot increase after blowup, 
we obtain that $h^0(mK^{-1})=0$ for $\tilde S_t$ in the case $n>4$ and $k\not\in\{1,n-1\}$.
Hence $h^0(mK^{-1})=0$ also for the surface $S_t$. 
\endproof

We recall that in Proposition \ref{prop:sub1}
we have obtained $G(k,n-k)$-equivariant deformation
of the LeBrun orbifold.
Correspondingly we have  $G(k,n-k)$-equivariant
deformation of the twistor space $Z_{\LB}$.
For algebraic dimension of these twistor spaces,
by using Propositions \ref{prop:dfmpair} and \ref{prop:dfmS} we obtain the following

\begin{proposition}\label{prop:ad}
(i) If $k\in\{1,n-1\}$, all the deformed twistor spaces are Moishezon.
at least for small deformations.
(ii) If $n=4$ and $k=2$, there exists a
small deformations whose  algebraic dimension 
is two.
(iii) If $n>4$ and $k\not\in\{1,n-1\}$, 
there exists a small deformation whose algebraic dimension is zero.
\end{proposition}

\proof
%
%
For (i), from the concrete construction of the equivariant deformations
of the surface $S$, we readily see that 
if $k\in \{1,n-1\}$ the deformed surface $S_t$ has 
a unique U(1)-invariant real anticanonical curve.
Further the curve is a cycle of smooth
rational curves consisting of four
irreducible components, regardless of 
the value of $n$.
Let $Z_t$ be the twistor space in which the surface $S_t$ is contained.
The exact sequence
$0\to \ms O\to K^{-1/2}\to K^{-1}_{S_t}\to 0$ on $Z_t$ means that the system $|K^{-1/2}|$ contains a U(1)-invariant pencil whose base curve is exactly the above cycle.
(When $n>3$ the pencil is exactly $|K^{-1/2}|$ itself.)
Further as any element of this pencil contains 
the cycle, general members of the pencil
also satisfy $\kappa^{-1}=2$.
Thus the twistor space $Z_t$ has a pencil
whose general members satisfy $\kappa^{-1}=2$.
This directly implies that $Z_t$ is Moishezon
\cite{U}.
Further, the twistor spaces
 obtained as the
$G(k,n-k)$-equivariant deformation of the pair $(Z_{\LB},S)$ exhausts the 1-dimensional family obtained in Proposition \ref{prop:sub1}
at least for small deformations,
 since the deformation of the pair actually deforms the complex structure of $Z_{\LB}$
by Proposition \ref{prop:dfmpair}.
Hence we obtain the assertion (i).
The assertion (ii) can be obtained in a similar way.
(iii) is much easier.
\proofend
%

\begin{remark}
{\em
In (iii) of Proposition \ref{prop:ad}
it is very likely that all the twistor spaces
in the 3-dimensional family have algebraic dimension zero, at least for small deformations.
}
\end{remark}

Now we can give a proof of 
Theorem \ref{thm:sfk} in the introduction.

\vsp
\noindent 
{\em Proof of Theorem \ref{thm:sfk}.}
As in the proof of Theorem \ref{thm:smry1} by restricting the 
versal family 
for the twistor space $Z_{\LB}$ on 
$\widehat{\ms O(-n)}$
to  the 
U(2)-invariant subset $B_1$ 
in $S^{n-2}_1\CC^2\oplus S^{n-4}_2\CC^2$,
we obtain the versal family of 
$G(n-1,1)$-equivariant deformations of 
$Z_{\LB}$.
Take any  one-dimensional subspace of 
$\langle z^{n-2}\rangle_{\CC}$ over $\RR$
and consider the restriction of the 
versal family to this subspace.
We show that
on 
the open subset $\ms O(-n)$
 the corresponding family of
ASD structures  provides
the required family of ASD ALE K\"ahler metrics.

Let $Z_t$ be any one of the twistor spaces
in this real one-dimensional family,
and $L_{\infty}\subset Z_t$ the twistor line over the 
orbifold point.
We need to show that the linear system 
$|K^{-1/2}|$ on $Z_t$ carries a real reducible member which contains $L_{\infty}$.
Recall from the proof of Proposition \ref{prop:ad} (i) 
that the linear system $|K^{-1/2}|$ on $Z_t$  has a U(1)-invariant  
pencil, and the base locus
of the pencil is a cycle of four rational curves.
Take a uniformization of the orbifold point 
$\infty\in\widehat{\ms O(-n)}$, and let $\tilde{\infty}$ be the point over
$\infty$, and  $\Gamma\simeq\ZZ_n$  the group at
$\tilde{\infty}$.
Take a $\Gamma$-invariant open neighborhood $V$ of the twistor line
$L_{\infty}$, and let $u:\tilde V\to V$ be the uniformization corresponding
to the above uniformization on the base.
Let $ L_{\tilde\infty}=u^{-1}(L_{\infty})$ be the 
$\Gamma$-invariant twistor line
over  $\tilde{\infty}$.
Then in the present situation 
the group $\Gamma$ acts on the normal bundle $N=N_{L_{\tilde\infty}/\tilde V}$ 
as 
\begin{align}\label{act10}
(x, y, z)\longmapsto (\zeta x,\zeta^{-1} y,z)
\end{align} where $(x,y)$ is  holomorphic fiber coordinates on $N$, 
$z$ is a coordinate on $L_{\tilde\infty}$, and $\zeta=e^{2\pi i/n}$ is the generator of $\Gamma
$ as before.
We consider the pullback of the
U(1)-invariant pencil in $|K^{-1/2}|$ on $Z_t$ to $\tilde V$
by the uniformization map $u$.
Since $u^*K^{-1/2}_{\tilde V}\simeq K^{-1/2}_V$, 
the pullback is a pencil whose members belong to $|K^{-1/2}|$ 
of the open twistor space $\tilde V$.

Now adapting the argument of Kreussler given in the proof of \cite[Proposition 3.7]{Kr98},
we consider 
an element $S_0$ of the last pencil on $\tilde V$, 
which is uniquely specified by the property that 
it goes through a generic point of $L_{\tilde\infty}$.
Here, genericity means that the point does not belong
to the base curve of the pencil.
Then by that argument, the unique divisor $S_0$ is of the form $D_0+\ol D_0$, 
where $D_0$ and $\ol D_0$ are irreducible 
non-singular, and intersect along 
$ L_{\tilde\infty}$ transversally.  
(Note that Pedersen-Poon's result about reducibility of certain divisor 
used in the Kreussler's argument 
does not require for the twistor space
to be compact.)
Taking the image of $S_0$ under the uniformization map $u$, it follows
that a member of the
U(1)-invariant pencil on the original twistor space $Z$
which goes through the  generic point of $L_{\infty}$ is unique, and it
contains the whole of $L_{\infty}$.
Moreover, the last member is clearly real,
and at least on the neighborhood $V$,
it consists of two irreducible components
$u(D_0)$ and $u(\ol D_0)$
whose intersection is precisely $L_{\infty}$.

We now show that this divisor is reducible on the whole of $Z_t$. 
The  $\Gamma$-action on the uniformization $\tilde V$ preserves
each of the two irreducible components of the divisor $S_0$, and from 
\eqref{act10}, it acts on each of the components as merely as a  multiplication by $\zeta$
in the normal direction. 
This means that the images 
$u(D_0)$ and $u(\ol D_0)$
are non-singular, and the self-intersection numbers of $L_{\infty}$ in
the components are both $(+n)$.
Therefore, if the divisor is irreducible, its normalization would have two disjoint 
curves whose self-intersection numbers are both $(+n)$.
This contradicts Hodge index theorem.
Therefore the divisor is reducible on $Z_t$.
Hence by the theorem of Pontecorvo, we obtain that on the smooth locus
$\ms O(-n)$, the ASD structure
associated to the twistor space $Z_t$ is represented by a K\"ahler metric.
Also the presence of the above divisor $D_0 + \ol D_0$, as well as the above 
degree of the normal bundle of $L_{\infty}$ in $D_0$ and $\ol D_0$,
mean that the K\"ahler metric is ALE at infinity
(see \cite[Proposition 6]{LB92}).

Finally we detect the complex structure on 
the regular locus $\ms O(-n)$. 
For this let $D_t+\ol D_t$ be the reducible member of $|K^{-1/2}|$ on $Z_t$
obtained above.
Then as both $D+\ol D$ and $D_t+\ol D_t$ are
unique reducible member containing the twistor
line over the orbifold point $\infty$,
the divisor
 $D_t+\ol D_t$ is naturally regarded as a deformation of the 
divisor $D+\ol D$ on $Z_{\LB}$.
We determine the complex structure of $D_t$.
For this, since $D_t+\ol D_t\in |K^{-1/2}|$, we have, by adjunction formula,
$
K_{D_t} = K_{Z_t} + D_t|_{D_t} = K^{1/2} -\ol D_t|_{D_t}.
$
Hence, since $\ol D_t|_{D_t}\simeq\ms O_{D_t}(L_{\infty})$, we have 
\begin{align}\label{adj1}
K_{D_t}^{-1}\simeq K^{-1/2}|_{D_t} + L_{\infty}.
\end{align}
Let $S_t\in |K^{-1/2}|$ be any member
of the U(1)-invariant pencil which is different from $D_t+\ol D_t$.
Then we can write $S_t\cap D_t=C_1+C_2$ for two components $C_1$ and $C_2$ of 
the base locus of the pencil.
By \eqref{adj1} the curve $C_1+C_2+L_{\infty}$ is an anticanonical curve
on $D_t$, and it is a triangle.
On the other hand, noting that on $\FF_n$
the $(-n)$-section is a base curve of the anticanonical system (as $n>2$) and it is disjoint from any $(+n)$-section,
$\FF_n$ does not have such a triangle anticanonical curve. 
Therefore $D_t$ is not biholomorphic to $\FF_n$.

Now we show from these that $D_t$ $(t\neq 0)$ is biholomorphic to $\FF_{n-2}$.
For this we recall that any small deformation of  rational ruled surface $\FF_n$
must be of the form $\FF_{n-2k}$ where $k\ge 0$ and $n-2k\ge 0$
(see \cite{Suwa} for the Kuranishi family of $\FF_n$.)
Also, in these deformations the ruling (i.e.\,the projection to $\CP^1$) is preserved.
Then as the pair $(D_t,L_{\infty})$ is obtained 
as a small deformation of the pair $(D_0,L_{\infty})$ which satisfies $L_{\infty}^2 = n$,  
we still have  $L_{\infty}^2 = n$ on $D_t$.
This readily means that on $D_t\simeq \FF_{n-2k}$ we have 
\begin{align}\label{Linfty}
L_{\infty}\sim \Gamma +  (n-k) h \quad{\text{(linear equivalence)}}
\end{align}
where $h$ denotes the fiber class of the ruling, and $\Gamma$ denotes a section of the ruling that satisfies $\Gamma^2 = -(n-2k)$.
(Of course such a section is  unique as long as $n-2k>0$.)
On the other hand, on $\FF_{n-2k}$ we have $K^{-1}\sim 2\Gamma + (2+n-2k)h$. 
As $C_1 + C_2 + L_{\infty}\sim K^{-1}$ as above, it follows that we may suppose that 
\begin{align}\label{C1C2}
C_1\sim \Gamma + (1-k)h, \quad C_2\sim h.
\end{align}
But since $C_1$ is an irreducible curve, we have $1-k\ge 0$.
Hence as $k>0$ (since $D_t\not\simeq\FF_n$ as above), we obtain $k=1$.
Thus we have $D_t\simeq\FF_{n-2}$ and $L_{\infty}\in|\Gamma+(n-1)h|$.
It is not difficult to show that 
this linear system induces
an embedding $\FF_{n-2}\subset \CP^{n+1}$
whose image is a non-singular surface of degree $n$. Therefore the complement
$\FF_{n-2}\backslash L_{\infty}$ is an 
algebraic surface in $\CC^{n+1}$.
\endproof

\vsp
Theorem \ref{thm:sfk} gives a partial answer to a question by 
Viaclovsky \cite[1.4 Question (2)]{V}
concerning scalar-flat K\"ahler deformations of the LeBrun's ALE metric.
Moreover, in relation with 
 SU(2)-invariant scalar-flat K\"ahler metrics
obtained by Pedersen-Poon \cite{PP90},
Theorem \ref{thm:sfk} shows that their metrics are incomplete
as long as they are obtained as small deformations of
the LeBrun's ALE metric.
On the other hand, when $k\not\in\{1,n-1\}$,
the $G(k,n-k)$-equivariant deformation
obtained in Proposition \ref{prop:sub1}
does not preserve K\"ahlerity in general,
since under the $G(k,n-k)$-equivariant deformation
of the pair $(Z_{\LB}, S)$ obtained in 
Proposition \ref{prop:dfmpair},
the divisor $D+\ol D$ can be shown to disappear
from the structure of the deformed surface.

Because the twistor spaces obtained in 
Theorem \ref{thm:sfk} 
possess the irreducible singular member
in $|K^{-1/2}|$ whose minimal resolution is exactly the one appeared in 
\cite{Hon07-2,HonDSn1}, it is very natural to 
expect that these twistor spaces
on $\widehat{\ms O(-n)}$ has a structure of
a double cover of $\CP^3$ (if $n=3$)
or a scroll of planes in $\CP^n$ (if $n>3$), 
whose branch divisor 
is a quartic surface (if $n=3$) or 
a cut of the scroll by a quartic hypersurface (if $n>3$).
In other words, it is quite likely that 
the ASD K\"ahler metrics on $\ms O(-n)$ in Theorem \ref{thm:sfk}
could be obtained as a limit of the ASD structures
on $n\CP^2$ which correspond to the twistor spaces
investigated in \cite{Hon07-2,HonDSn1}.


\end{document}